\RequirePackage[l2tabu, orthodox]{nag}
\RequirePackage{snapshot}

\documentclass[10pt,onecolumn]{extarticle}

\makeatletter
\if@twocolumn
  \usepackage[dvips,letterpaper,top=0.5in, bottom=0.5in, left=0.75in, right=0.5in,includefoot,heightrounded]{geometry}
\else
  \usepackage[dvips,letterpaper,margin=1in,includefoot,heightrounded]{geometry}
\fi

\usepackage{srcltx}

\usepackage[russian,portuges,english]{babel}

\iflanguage{portuges}
    {\newcommand{\keywordname}{Palavras-chaves}}
    {\newcommand{\keywordname}{Keywords}}

\usepackage{amsmath}
\usepackage{amssymb,amsfonts}

\usepackage{abstract}

\usepackage{graphicx}
\usepackage{epsfig}
\usepackage[usenames,dvipsnames,svgnames,x11names]{xcolor}
\usepackage{subfigure}

\usepackage{booktabs}

\usepackage{setspace}
\usepackage{flushend}

\usepackage{cite}

\usepackage{hyperref}
\usepackage[normalem]{ulem}
\usepackage{bm}

\usepackage{enumerate}

\usepackage[noend]{algpseudocode}

\usepackage{listings}

\usepackage[short,12hr]{datetime}

\usepackage[sort&compress,square,comma,authoryear]{natbib}

       \usepackage{fouriernc}
\makeatletter

\newcommand{\printtitle}{%
\makeatletter
\if@twocolumn

\twocolumn[%
  \maketitle
  \begin{onecolabstract}
    \myabstract
  \end{onecolabstract}
  \begin{center}
    \small
    \textbf{\keywordname}
    \\\medskip
    \mykeywords
  \end{center}
  \bigskip
]
\saythanks
\else
  \maketitle
  \begin{onecolabstract}
    \myabstract
  \end{onecolabstract}
  \begin{center}
    \small
    \textbf{\keywordname}
    \\\medskip
    \mykeywords
  \end{center}
  \bigskip
  \onehalfspacing
\fi
\makeatother
}

\author{%
J. Le\~ao%
\thanks{Departamento de Estat\'istica, Universidade Federal do Piau\'i, Teresina, Brazil.}
\and
H. Saulo%
\thanks{Departamento de Economia, Universidade Federal do Rio Grande do Sul, Porto Alegre, Brazil.}
\and
M. Bourguignon%
\thanks{Programa de Pós-gradua\c{c}\~ao em Estat\'istica, Universidade Federal de Pernambuco, Recife, Brazil.}
\and
R. J. Cintra%
\thanks{Departamento de Estat\'istica, Universidade Federal de Pernambuco, Recife, Brazil. E-mail: \url{rjdsc@de.ufpe.br}}
\and
L. C. R\^ego%
\thanks{Formely Departamento de Estat\'istica, Universidade Federal de Pernambuco, Recife, Brazil; currently
Departamento de Estat\'{\i}stica e Matem\'atica Aplicada,
Universidade Federal do Cear\'a, Brazil.
E-mail: \url{leandro@dema.ufc.br}}
\and
G. M. Cordeiro%
\thanks{Departamento de Estat\'istica, Universidade Federal de Pernambuco, Recife, Brazil. E-mail: \url{gauss@de.ufpe.br}}
}

\title{%
On Some Properties of the Beta Inverse Rayleigh Distribution}

\newcommand{\myabstract}{%
We study with some details a lifetime model of the class of beta generalized models, called the beta
inverse Rayleigh distribution, which is a special case of the Beta Fr\'echet distribution.
We provide a better foundation for some properties including quantile function, moments, mean deviations,
Bonferroni and Lorenz curves, R\'enyi and Shannon entropies and order statistics. We fit the proposed model
using maximum likelihood estimation to a real data set to illustrate its flexibility and potentiality.
}

\newcommand{\mykeywords}{%
Beta-Generated class,
Entropy,
Generalized distribution,
Maximum likelihood estimation,
Moment.
}

\date{}

\begin{document}

\printtitle

\section{Introduction}

After its inception by \cite{t:64}, the inverse Rayleigh (IR) distribution was championed
by \cite{vd:72} and \cite{il:73} during the 1970s. In \cite{vd:72}
several of its statistical properties were addressed, in particular, maximum likelihood estimation,
confidence intervals, and hypotheses tests.
An early application involved lifetime modeling of experimental units. More recently, \cite{g:93}
provided closed-form expressions for the mean, harmonic mean, geometric mean, mode and the median of this distribution.
In \cite{mhs:05} the negative moment estimator for the IR distribution was investigated.
Moreover, different methods of estimation have been numerically compared in \cite{g:93}
and \cite{sl:10}.
Acceptance sampling techniques also received a treatment based on the IR distribution \cite{rosh:05}.
In 2010, a model for lower record value based on the IR distribution was proposed in \cite{sl:10}
and a Bayesian approach for its associate parameter estimation.

Distribution generalization theory has received considerable attention in the past decades, \cite{Am:25},
\cite{Gd:53}, \cite{HW:87} and \cite{McD:84}.
A particular prominent generalization model is the class of beta generalized distributions, first
introduced in \cite{Eug:02}.  In this seminal work, the authors introduced the new class of distributions
from the logit of the beta random variable, and obtained as a special case the beta normal (BN)
distribuiton. This distribution could provide flexible shapes including bimodality,
being therefore a candidate for a wide range of applications.
Additional properties of the BN distribution have been studied in detail
by \cite{GN:04} and \cite{rcG:12}.
In a similar manner, other beta generalizations have been proposed taking into account several
baseline distributions. To cite a few, we identify
the beta Gumbel by \cite{NKt:04},
beta Fr\'echet by \cite{NGpt:04},
beta exponential by \cite{NKtz:05},
beta Weibull by \cite{LT:07},
beta Pareto by \cite{Ak:08},
beta generalized exponential by \cite{BSt:10},
beta generalized normal by \cite{CTa:11}
and beta generalized half-normal by \cite{Pta:10} distributions.

In this paper, we study the beta generalized distribution based on the IR distribution,
called the beta inverse Rayleigh (BIR) distribution. The BIR distribution is a special case of the
beta Fr\'echet (BF) distribution, which was introduced by \cite{NGpt:04} and studied by \cite{BSt:11}.
These two papers provide some mathematical properties for the BF distribuion, which in turn can be easily
adapted for the BIR distribution. We provide a better foundation for these and other mathematical
properties. An application to a real life data set is presented. The BIR
distribution is expected to have immediate application in reliability and survival studies.

The rest of the paper unfolds as follows. In Section 2,
we present the BIR distribution, derive its density and some expressions
for the cumulative distribution function~(cdf), and provide an analytical
study of the unimodality region. In Section 3,
we give the hazard rate function and its asymptotic behavior.
In Section 4, we derive the formulae for the moments.
Further, in Sections 5--8,
we derive quantile function, skewness and kurtosis, mean deviations, R\'enyi
entropy, Shannon entropy and order statistics.
In Section 9, we
discuss maximum likelihood estimation and
present the elements of the observed information matrix.
An application to real data is performed in Section 10.
Finally, in Section 11,
we offer some concluding remarks.

\section{The BIR distribution}

Let $G(x)$ be a baseline cumulative distribution function (cdf). Then, the associated beta generalized distribution
$F(x)$ based on the logit of the beta random variable is given by \citep{Eug:02}
\begin{equation}
\label{eq1}
F(x)=\operatorname{I}_{G(x)}(a,b),
\end{equation}
where $a>0$, $b>0$, $\operatorname{I}_{y}(a,b)$ is the incomplete beta function ratio
\begin{align*}
\operatorname{I}_{y}(a,b) = \frac{1}{\operatorname{B}(a,b)}
\int^{y}_{0} \omega^{a-1}(1-\omega)^{b-1}\mathrm{d}\omega,
\end{align*}
and $\operatorname{B}(\cdot,\cdot)$ denotes the beta function.
The extra shape parameters $a$ and $b$ control skewness, kurtosis and tail weights.

\begin{figure}
\centering
\subfigure[\label{Figure1a}]{\epsfig{file=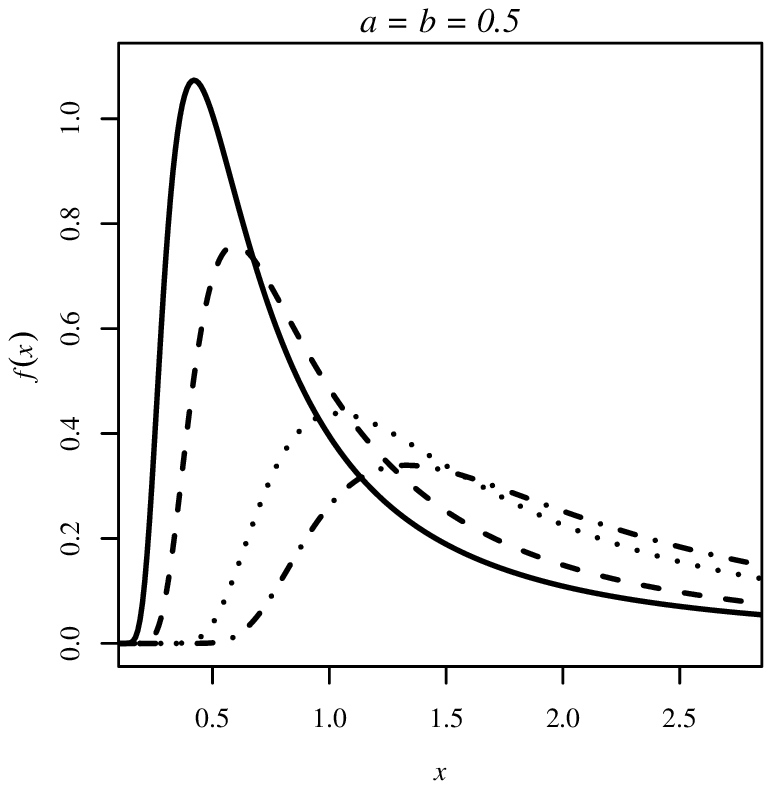,height=0.45\linewidth,angle=-90}}
\subfigure[\label{Figure1b}]{\epsfig{file=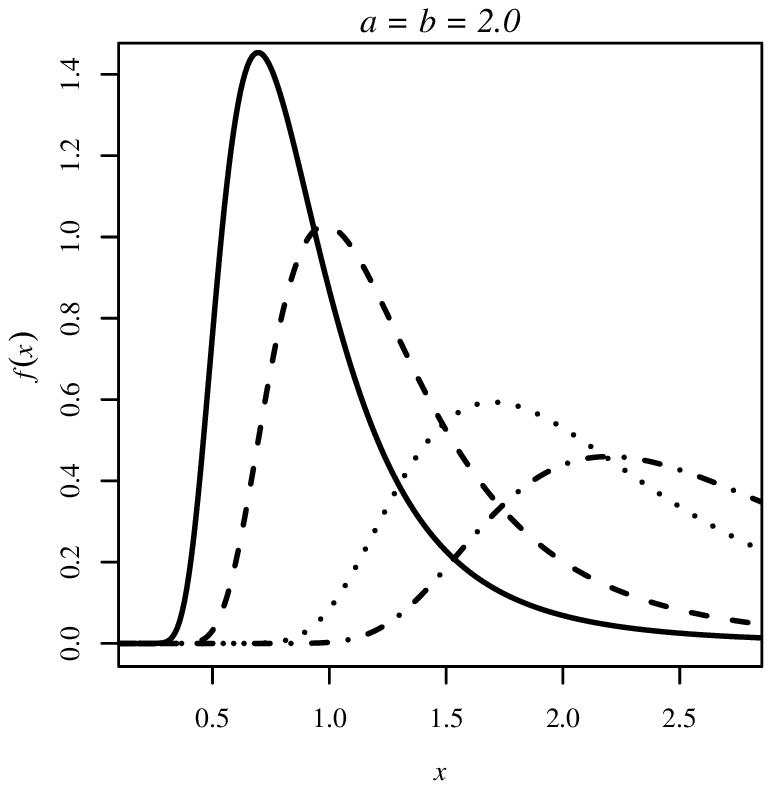,height=0.45\linewidth,angle=-90}}
\subfigure[\label{Figure1c}]{\epsfig{file=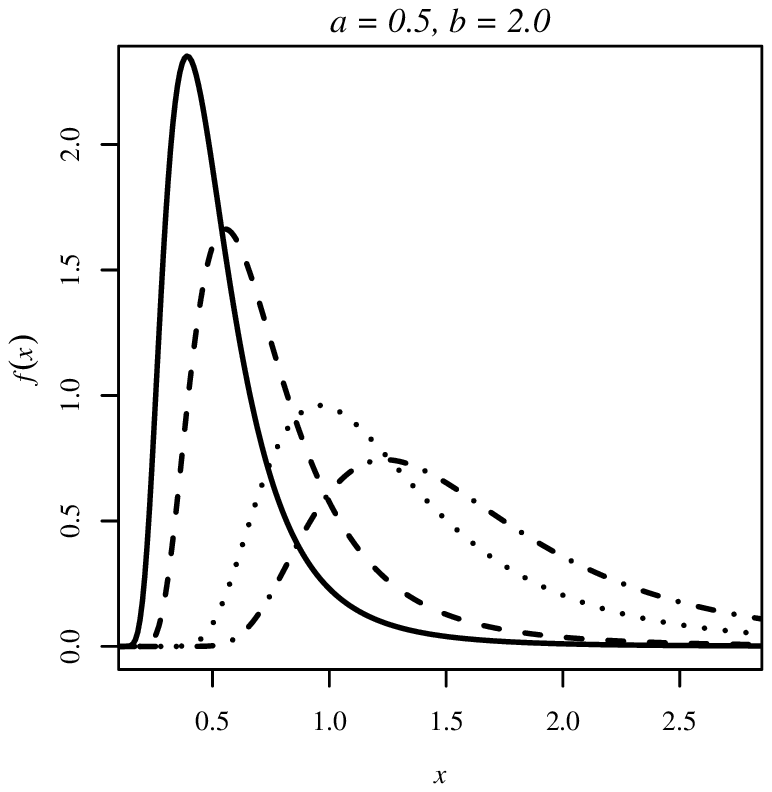,height=0.45\linewidth,angle=-90}}
\subfigure[\label{Figure1d}]{\epsfig{file=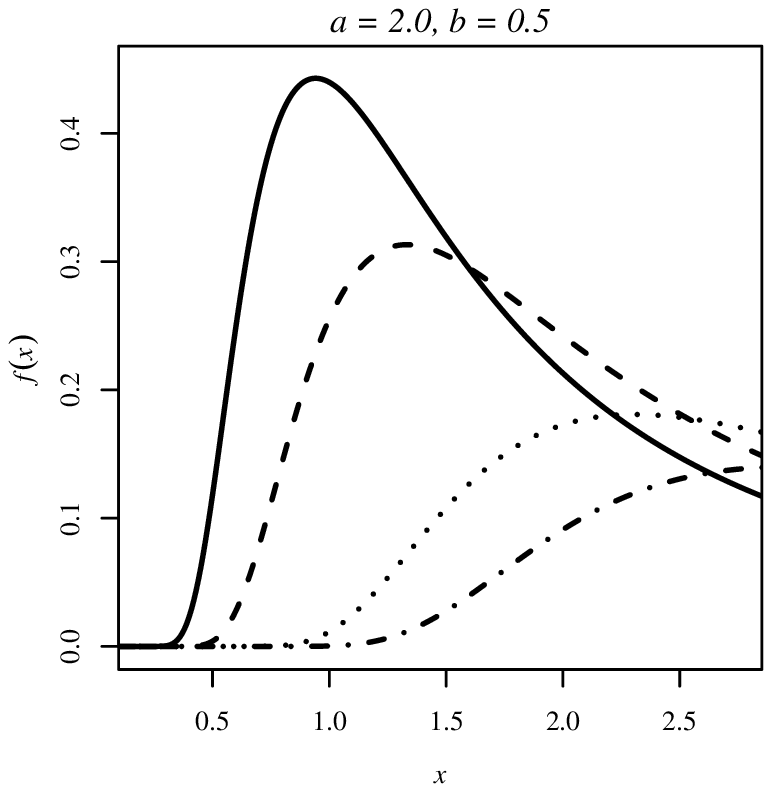,height=0.45\linewidth,angle=-90}}
\hfill\caption{Plots of the BIR pdf for $\theta = 0.5$ (solid line),
$\theta = 1.0$ (dashed line), $\theta = 3.0$ (dotted line) and
$\theta = 5.0$ (bold line).}
\label{Figure1}
\end{figure}

The IR distribution is a single-parameter distribution defined over the semi-infinite interval $[0,\infty)$.
Its cdf is given by
\begin{align*}
G(x;\theta) = \exp\left(-\frac{\theta}{x^2}\right), \quad x> 0,
\theta>0 .
\end{align*}
Inserting $G(x;\theta)$ into~\eqref{eq1},
we obtain the BIR cumulative distribution
\begin{align}
\label{eq2}
F(x) =
\operatorname{I}_{\exp\left(-\frac{\theta}{x^2}\right)}(a,b)
=
\frac{1}{\operatorname{B}(a,b)}
\int_0^{\exp\left(-\frac{\theta}{x^2}\right)}
\omega^{a-1}(1-\omega)^{b-1}\mathrm{d}\omega
,
\end{align}
for $x>0$, $a>0$, $b>0$ and $\theta>0$. Note that if we take Fr\'echet cdf
$G(x,\sigma,\lambda)=\exp\left\{-\left(\frac{\sigma}{x}\right)^{\lambda}\right\}$, where
$\sigma>0$ and $\lambda>0$ are the scale and shape parameters, respectively, into~\eqref{eq1}, we
obtain the BF distribution. Thus, the BIR model is obtained for $\sigma^{2}=\theta$ and $\lambda=2$.
Note also that for the special case $a=b=1/2$, the BIR cumulative function
has a closed-form expression given by
\begin{align*}
F(x) = \frac{2}{\pi}\arcsin \left\{\exp\left(-\frac{\theta}{2 x^2}\right)\right\}.
\end{align*}
The BIR probability density function (pdf) can be expressed as (for $x>0$)
\begin{align}\label{eq3}
f(x)
=
\frac{2\theta}{\operatorname{B}(a,b)\,x^3}
\exp\left(-\frac{a\theta}{x^2}\right)
\left[
1-\exp\left(-\frac{\theta}{x^2}\right)
\right]^{b-1}.
\end{align}
The BIR random variable $X$ is denoted by $X\sim\operatorname{BIR}(a,b,\theta)$.
The parameters $a$ and $b$ affect the skewness of $X$ by changing the relative tail weights. Figure~\ref{Figure1}
displays the BIR pdf for several choices of parameter values. Simulating the BIR random variable is
relatively simple. Let $Y$ be a random variable distributed according to the usual beta distribution with
parameters $a$ and $b$. Thus, by means of the inverse transformation method,
the random variable $X$ given by
\begin{align*}
X= \sqrt{ -\frac{\theta}{\log (Y)}}
\end{align*}
follows (\ref{eq3}).

\subsection{General expansion}

Although the cdf and pdf of $X$ require mathematical functions that are widely available in contemporary statistical
packages, \cite{EaT:02} and \cite{RCT:11} often further analytical and numerical derivations take advantage of
power series expansions for the cdf. From the BIR density function~\eqref{eq3}, the cdf of $X$ can be expressed after usual integration as

\begin{align*}
F(x) =
\frac{1}{\operatorname{B}(a,b)}
\int^{x}_{0}
\frac{2\theta}{y^{3}}
\exp\left(-\frac{a\theta}{y^{2}}\right)
\left\{
1 - \exp\left(-\frac{\theta}{y^{2}}\right)
\right\}^{b-1}
\mathrm{d}y.
\end{align*}
Setting $u=\theta y^{-2}$, it follows that
\begin{align}
\label{eq6}
F(x) = \frac{1}{\operatorname{B}(a,b)}
\int^{\infty}_{\frac{\theta}{x^2}}
\exp(-au)\left\{ 1-\exp(-u) \right\}^{b-1}
\mathrm{d}u
.
\end{align}
Notice that for $|z|<1$ and $b>0$ a real non-integer number,
we have the power series expansion
\begin{align}
\label{eq7}
(1-z)^{b-1}
=
\sum^{\infty}_{n=0}
\frac{(-1)^{n}\,\Gamma(b) }{\Gamma(b-n)\,n!}
z^n
,
\end{align}
where $\Gamma(\cdot)$ is the gamma function.
Applying this identity into~\eqref{eq6} yields
\begin{align*}
F(x)
=&
\frac{1}{\operatorname{B}(a,b)}
\sum^{\infty}_{n=0}
\frac{(-1)^n\Gamma (b)}{\Gamma(b-n)n!}
\int^{\infty}_{\frac{\theta}{x^2}}
\exp\{-(a+n)u\}
\mathrm{d}u
\nonumber
\end{align*}
and then
\begin{align}
\label{cdf-expansion}
F(x)=
\frac{1}{\operatorname{B}(a,b)}
\sum^{\infty}_{n=0}
\frac{(-1)^n\,\Gamma(b)}
{(a+n)\,\Gamma(b-n)\,n!}
\exp\left\{-\frac{(a+n)\theta}{x^2}\right\}
.
\end{align}

Now, considering the following quantity,
\begin{align*}
c_n(a,b)
=
\frac{(-1)^n\Gamma(a+b)}{(a+n)\Gamma(a)\Gamma(b-n)n!}
,
\end{align*}
we can write the BIR cdf as a linear combination of IR cdfs.
Indeed, we obtain
\begin{align*}
F(x) = \sum_{n=0}^\infty c_n(a,b)
\, G(x;(a+n)\theta).
\end{align*}

In a similar way, the BIR pdf can be expressed according
to the following linear combination
\begin{align*}
f(x) = \sum_{n=0}^\infty
c_n(a,b)
\,
g(x;(a+n)\theta),
\end{align*}
where $g(x;(a+n)\theta)$ denotes the IR density function with
parameter $(a+n)\theta$.

\subsection{Unimodality}

\begin{figure}
\centering
\includegraphics[scale=0.9,angle=-90]{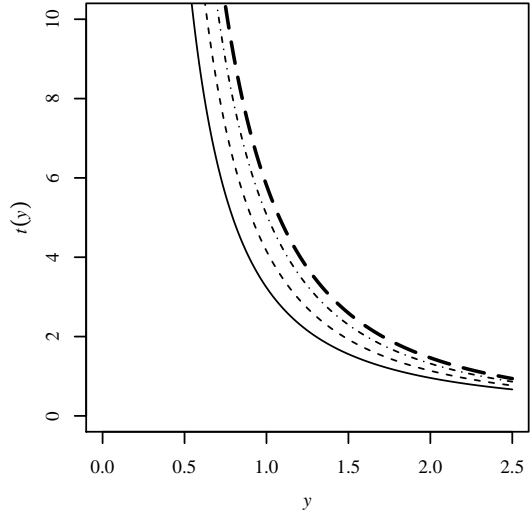}
\caption{Plots of $t(y)$ for $b = 0.25$ (solid line), $b = 0.5$
(dashed line), $b = 0.75$ (dotted line) and $b = 0.95$ (bold
line).}
\label{Figure2}
\end{figure}

The BIR distribution is unimodal for all values of $a,b,\theta > 0$.
In order to investigate the critical points of its density function,
the first derivative of $f(x)$ with respect to $x$ is given by

\begin{equation}
\label{derivative_of_f}
\begin{split}
\frac{d}{dx}
f(x) =
&
\frac{\theta^{2}}{B(a,b)\,x^6}
\exp\left(-\frac{a\theta}{x^2}\right)
\left[
1-\exp\left(-\frac{\theta}{x^2}\right)
\right]^{b-1}
\\
&
\times
\left[
4a
-
\frac{6x^2}{\theta}
+
\frac{4(b-1)}{1-\exp\left(\frac{\theta}{x^2}\right)}
\right]
,
\quad
x>0
.
\end{split}
\end{equation}
The signal of this derivative is determined by the expression in the last square brackets, since
the remaining terms are all positive. Considering the substitution $y=\theta x^{-2}$, the expression in
square brackets becomes
\begin{align}
\label{bracket_term}
4a
-
\frac{6}{y}
+
4\frac{b-1}{1-\exp(y)}
.
\end{align}
Now, we demonstrate that this expression is a monotonic function; therefore, \eqref{derivative_of_f} has
a single zero, which implies a unique mode.
Indeed, the derivative of~\eqref{bracket_term} becomes
\begin{align*}
t(y) = \frac{6}{y^2} + 4(b-1) \frac{\exp(y)}{\left[ 1-\exp(y)
\right]^2} .
\end{align*}
For $b\geq1$, this derivative is clearly positive.
For $0<b<1$,
Figure~\ref{Figure2} displays the numerical results that illustrate
the positiveness of the derivative of~\eqref{bracket_term}.

Moreover,
let $y_0$ be the zero of~\eqref{bracket_term}.
The BIR mode location is then given by
$\sqrt{\theta/y_0}$.
Since $y_0$ is independent of $\theta$,
the mode location is an increasing function of $\theta$.

\section{Hazard rate function}
\label{section-hazard}

The survival and hazard rate functions are given by $S(x)=1-F(x)$ and
$h(x)=f(x)/S(x)$, where $F(x)$ and $f(x)$ are the BIR cdf
and pdf, respectively. Thus, the hazard rate function of the random variable $X$ is
\begin{align*}
h(x)= \frac{2\theta}{\operatorname{B}(a,b)\,x^3} %
\frac{\left[1-\exp\left(-\frac{\theta}{x^2}\right)\right]^{b-1}\exp\left(-\frac{a\theta}{x^2}\right)}
{I_{1-\exp\left(-\frac{\theta}{x^2}\right)}(b,a)}
.
\end{align*}
Notice that we applied in~\eqref{eq2}
the symmetry property of the
incomplete beta function
$1-\operatorname{I}_x(a,b) = \operatorname{I}_{1-x}(b,a)$.

We now examine
the asymptotic behavior of $h(x)$ when $x\to\infty$ or $x\to0$.
First,
we prove that
$h(x) \sim 1/x$ as $x\to\infty$.
To establish this result, we verify that
$\lim_{x\to\infty}h(x)/x^{-1}$ is a constant.
Indeed, we have
\begin{align*}
\lim_{x\to\infty}
\frac{h(x)}{1/x}
=
&
\lim_{x\to\infty}
\frac{2\theta}{\operatorname{B}(a,b)}
\frac{1}{x^3}
\frac{\left[ 1-\exp\left(-\frac{\theta}{x^2}\right) \right]^{b-1} \exp\left(-\frac{a\theta}{x^2}\right)}{\operatorname{I}_{1-\exp(-\theta/x^2)}(b,a)}
x.
\end{align*}
Since
$\exp\left(-\frac{a\theta}{x^2}\right)\to1$ as $x\to\infty$,
we can write
\begin{align*}
\lim_{x\to\infty}
\frac{h(x)}{1/x}
=
&
\frac{2\theta}{\operatorname{B}(a,b)}
\lim_{x\to\infty}
\frac{\left[ 1-\exp\left(-\frac{\theta}{x^2}\right) \right]^{b-1}/x^2}{\operatorname{I}_{1-\exp(-\theta/x^2)}(b,a)}
.
\end{align*}
For any value of $b>0$,
the last expression gives rise to an inderteminate form.
Invoking L'H\^opital's rule and again considering that
$\exp\left(-\frac{a\theta}{x^2}\right)\to1$ as $x\to\infty$, we obtain
\begin{align*}
\lim_{x\to\infty}
\frac{h(x)}{1/x}
=
&
2\theta (b-1)
\lim_{x\to\infty}
\left[
\frac{1/x^2}{1-\exp\left(-\frac{\theta}{x^2}\right)}
\right]
-
2
.
\end{align*}
Applying the L'H\^opital rule again
we note that the above limit is well-defined
and is equal to $-2b$.

\begin{figure}
\centering
\subfigure[\label{Figure3a}]{\epsfig{file=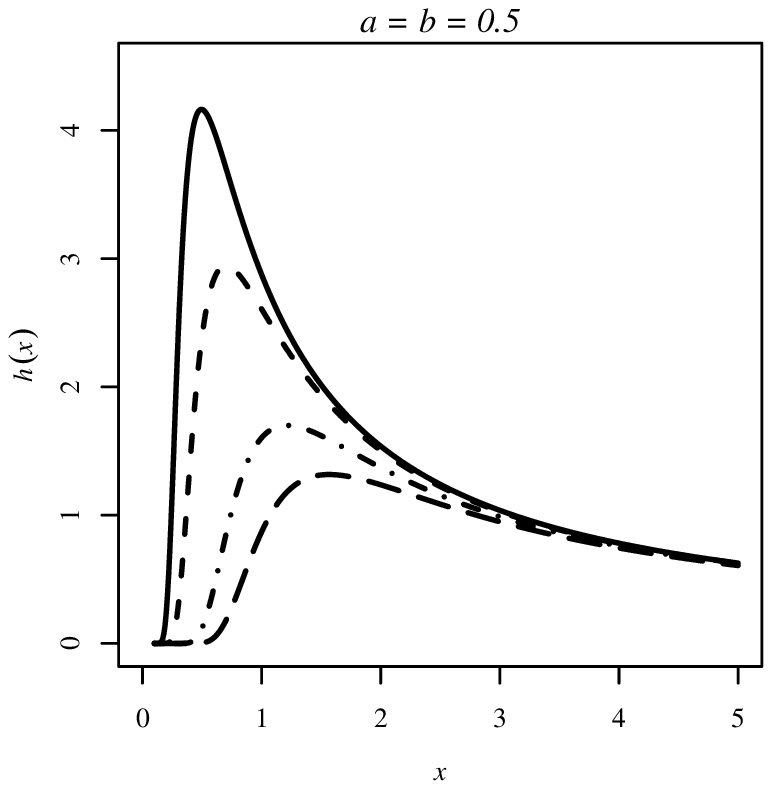,height=0.45\linewidth,angle=-90}}
\subfigure[\label{Figure3b}]{\epsfig{file=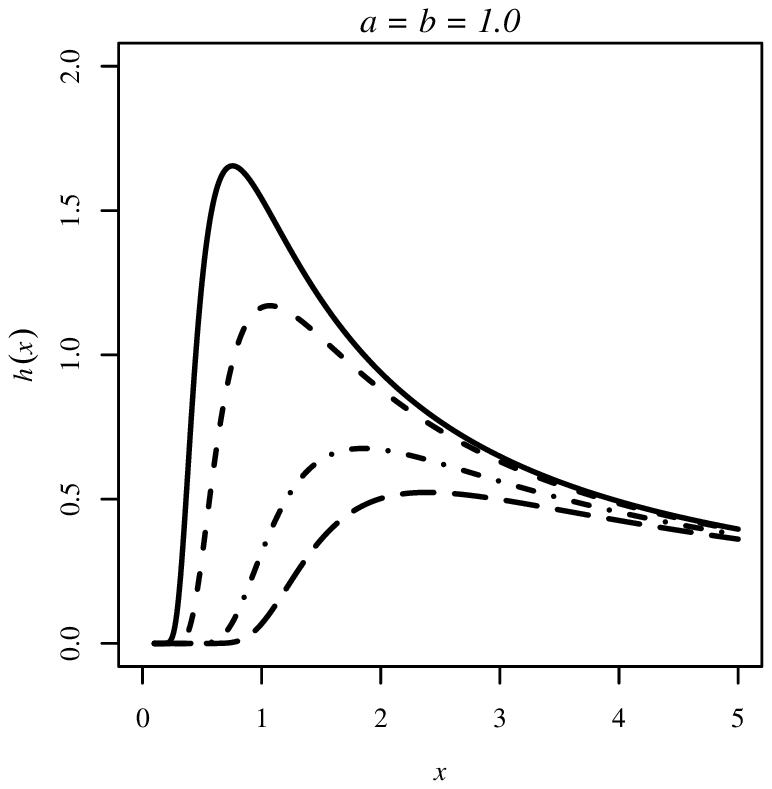,height=0.45\linewidth,angle=-90}}
\subfigure[\label{Figure3c}]{\epsfig{file=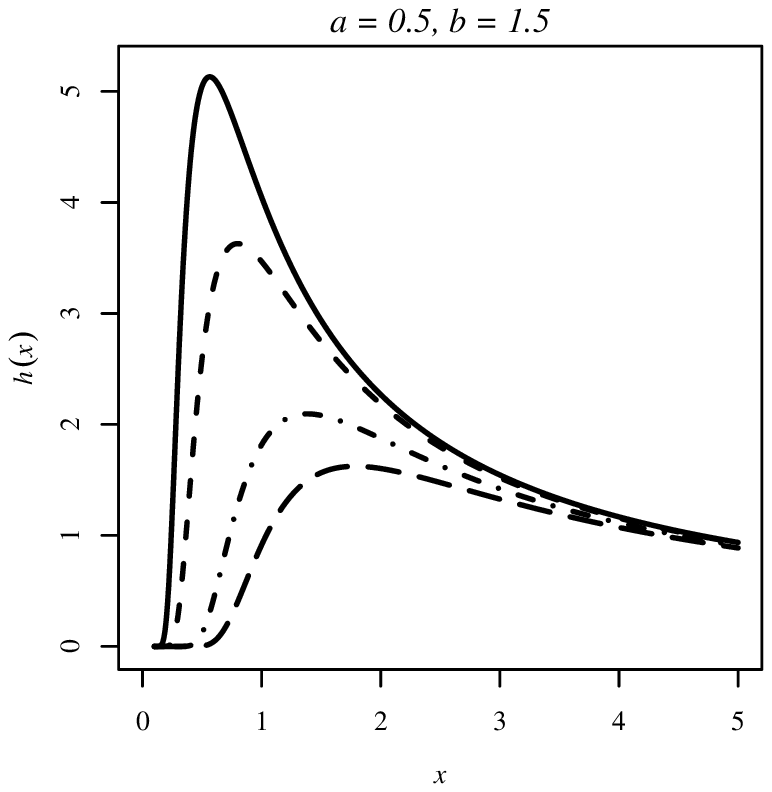,height=0.45\linewidth,angle=-90}}
\subfigure[\label{Figure3d}]{\epsfig{file=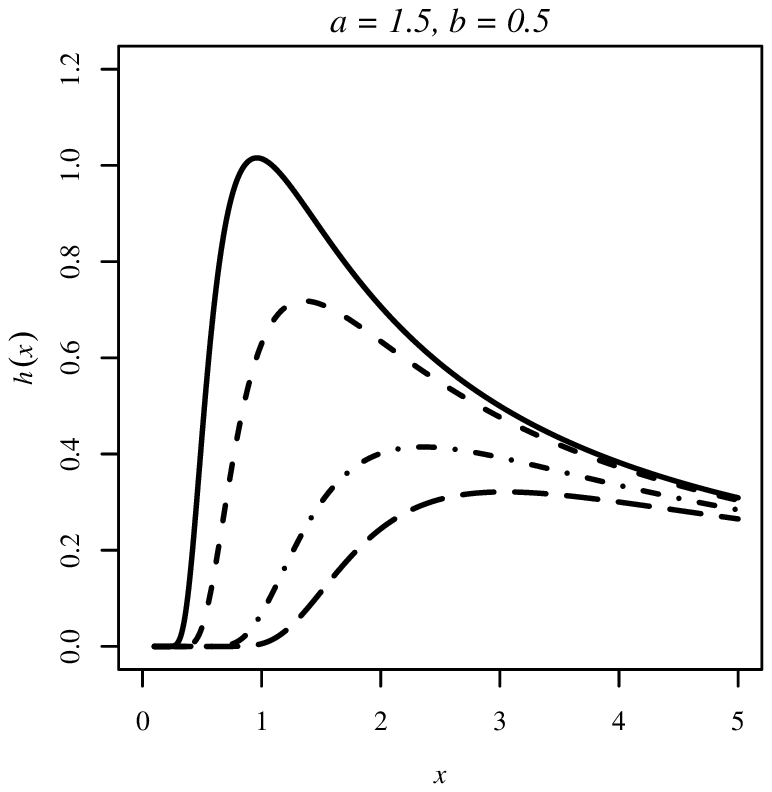,height=0.45\linewidth,angle=-90}}
\hfill\caption{Plots of the BIR hazard rate function
for
$\theta = 0.5$ (solid line),
$\theta = 1.0$ (dashed line),
$\theta = 1.5$ (dotted line)
and
$\theta = 3.0$ (bold line).}
\label{Figure3}
\end{figure}

Similarly,
let us show that
$h(x)\sim \exp(-a\theta/x^2) / x^3$
as $x\to0$.
In fact,
we have immedia\-tely that
\begin{align*}
\lim_{x\to0}
\frac{h(x)}{ \exp(-a\theta/x^2) / x^3}
=
&
\frac{2\theta}{\operatorname{B}(a,b)}
\lim_{x\to0}
\frac{\left[ 1-\exp\left(-\frac{\theta}{x^2}\right) \right]^{b-1}}{\operatorname{I}_{1-\exp(-\theta/x^2)}(b,a)}
\\
=
&
\frac{2\theta}{\operatorname{B}(a,b)}
.
\end{align*}
Notice also that
$\lim_{x\to0} \exp(-a\theta/x^2) / x^3 = 0$. Figure~\ref{Figure3} displays the behavior of $h(x)$
for selected values of the model parameters.

\section{Moments}

The moments play a crucial role in any statistical analysis.
The $r$th moment of $X$ is
\begin{align*}
\operatorname{E}(X^r)
=&
\frac{2\theta}{\operatorname{B}(a,b)}
\int_0^\infty
x^{r-3}
\exp\left(-\frac{a\theta}{x^2}\right)
\left[
1-\exp\left(-\frac{\theta}{x^2}\right)
\right]^{b-1}
\mathrm{d}x.
\end{align*}
Now,
we simplify the above integral.
First, letting $y=\theta x^{-2}$, we have
\begin{align*}
\operatorname{E}(X^r)
&=
\frac{\theta^{r/2}}{\operatorname{B}(a,b)}
\int_0^\infty
y^{-r/2}
\exp(-ay)
\left\{
1-\exp(-y)
\right\}^{b-1}
\mathrm{d}y
.
\end{align*}
We refer to the last integral
as $S_r(a,b)$.
Applying the series expansion~\eqref{eq7},
for any real $r$, we obtain
\begin{equation}
\label{eq.srab}
\begin{split}
S_r(a,b)
&=
\int_0^\infty
y^{-r/2}
\exp(-ay)
\sum_{n=0}^\infty
(-1)^n
\frac{\Gamma(b)}{\Gamma(b-n)n!}
\exp(-ny)
\mathrm{d}y
\\
&=
\sum_{n=0}^\infty
(-1)^n
\frac{\Gamma(b)}{\Gamma(b-n)n!}
\int_0^\infty
y^{-r/2}
\exp\left\{-(a+n)y\right\}
\mathrm{d}y
.
\end{split}
\end{equation}
This integral has a closed-form expression
by means of a direct application of
the gamma function integral, \citep{AbS:72}.
Since $a+n>0$, some manipulations yield
\begin{align}
\label{abramowitz611}
\int_0^\infty
y^{-r/2}
\exp\left\{-(a+n)y\right\}
\mathrm{d}y
=
\frac{\Gamma\left(1-\frac{r}{2}\right)}{(a+n)^{1-r/2}},
\quad
r<2.
\end{align}
Therefore, we can rewrite~\eqref{eq.srab} as
\begin{align*}
S_r(a,b)
=
\Gamma(b)
\Gamma\left(1-\frac{r}{2}\right)
\sum_{n=0}^\infty
\frac{(-1)^n}{(a+n)^{1-r/2}\Gamma(b-n) n!}
,
\quad
r<2
.
\end{align*}
If $b>0$ is an integer,
we obtain
\begin{align*}
S_r(a,b)
=
\Gamma\left(1-\frac{r}{2}\right)
\sum_{n=0}^b (-1)^n
\binom{b-1}{n} \frac{1}{(a+n)^{1-r/2}}
,
\quad
r<2
.
\end{align*}
We can write the $r$th moment of $X$ as
\begin{align*}
\operatorname{E}(X^r)
=
\frac{\theta^{r/2}}{\operatorname{B}(a,b)}
\,
S_r(a,b)
,
\quad
r<2.
\end{align*}
In particular,
for $r=1$ and integer an $b$,
we obtain
\begin{align*}
\operatorname{E}(X)
=
\frac{\sqrt{\pi\theta}}{\operatorname{B}(a,b)}
\sum_{n=0}^b \binom{b-1}{n} \frac{1}{\sqrt{a+n}}
.
\end{align*}
Negative moments can also be evaluated.
For example,
considering $r=-1$ and for an integer $b$,
we have
\begin{align*}
\operatorname{E}(X^{-1})
=
\frac{\sqrt{\pi/\theta}}{2\operatorname{B}(a,b)}
\sum_{n=0}^b (-1)^n
\binom{b-1}{n}
\frac{1}{\sqrt{(a+n)^3}}
.
\end{align*}

Notice that attempting to compute~\eqref{abramowitz611} outside $r<2$
gives undefined forms. For ins\-tance, if $r=2$, we have
\begin{align*}
\int_0^\infty
\frac{\exp[-(a+n)y]}{y}
\mathrm{d}y
=
\operatorname{E}_1(0),
\end{align*}
where $\operatorname{E}_1(\cdot)$ is the exponential integral
function \citep{AbS:72},
which tends to $-\infty$ as its argument goes to zero.
As a consequence,
the second moment of $X$ does not exist,
as well as all remaining higher order moments.
It is known that the second and higher order moments of
IR distribution are inexistent \citep{vd:72}.
As shown above,
the BIR distribution inherits this characteristic.

\section{Quantile function and quantile measures}
\label{section-quantile}

The quantile function of $X$ is given by
\begin{align*}
Q(u) = F^{-1}(u)
=
\sqrt{-\frac{\theta}{\log (\operatorname{I}^{-1}_{u}(a,b))}},
\quad 0<u<1,
\end{align*}
where $\operatorname{I}^{-1}_{u}(a,b)$ is
the inverse of the incomplete beta function.
The function $\operatorname{I}^{-1}_{u}(a,b)$ can be written as a
power series expansion \cite{WfA:11}
\begin{align*}
\operatorname{I}^{-1}_{u}(a,b)
=
\sum_{i=1}^{\infty}
q_i
\,
[a\operatorname{B}(a,b) u]^{i/a}
,
\end{align*}
where $q_1=1$ and
the remaining coefficients satisfy the following recursion
\begin{align*}
q_{i}
=&
\frac{1}{i^{2} + (a-2)i + (1-a)}
\left\{
(1-\delta_{i,2})
\sum_{r=2}^{i-1}
q_{r}q_{i+1-r}
[r(1-a)(i-r)-r(r-1)]
+
\right.
\\
&
\left.
\sum_{r=1}^{i-1}
\sum_{s=1}^{i-r}
q_{r}q_{s}q_{i+1-r-s}
[r(r-a)+s(a+b-2)(i+1-r-s)]
\right\}
,
\end{align*}
where $\delta_{i,2}=1$ if $i=2$ and $\delta_{i,2}=0$ if $i\neq2$.

Because the second, third, and fourth moments of the
BIR distribution are nonexistent,
usual skewness and kurtosis are not defined.
However,
quantile based measures, such as Bowley skewness \citep{kkP:62}
and Moors kurtosis \citep{MU:98},
can quantify asymmetry and the peakedness of a given distribution.
These measures exist even when moments are not available.
Bowley skewness and Moors kurtosis are expressed according to
\begin{align*}
B=&\frac{Q(3/4)-2Q(1/2)+Q(1/4)}{Q(3/4)-Q(1/4)},
\\
M=&\frac{Q(7/8)-Q(5/8)-Q(3/8)+Q(1/8)}{Q(6/8)-Q(2/8)}
.
\end{align*}

\begin{figure}
\centering
\subfigure[\label{Figure4a}]{\epsfig{file=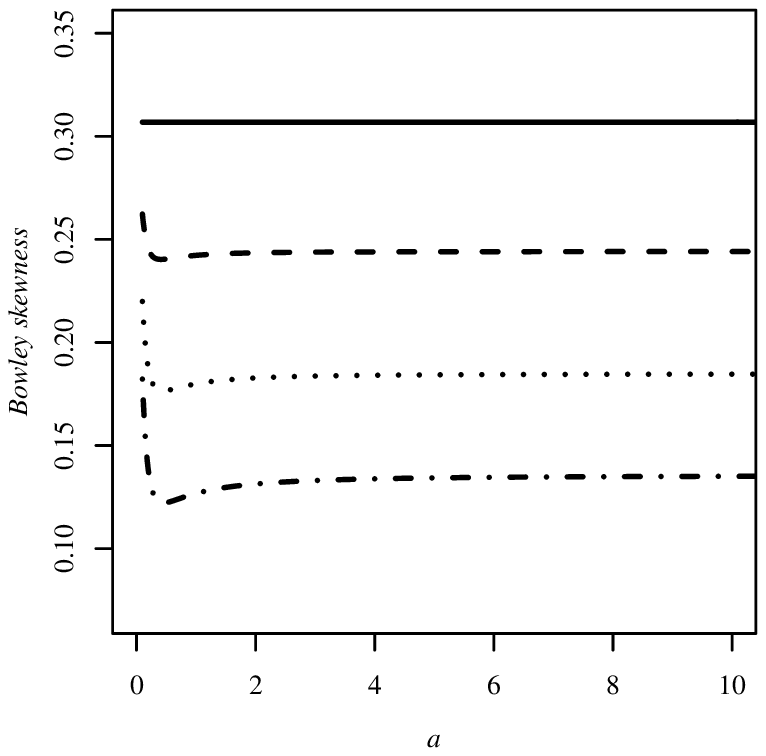,height=0.45\linewidth,angle=-90}}
\subfigure[\label{Figure4b}]{\epsfig{file=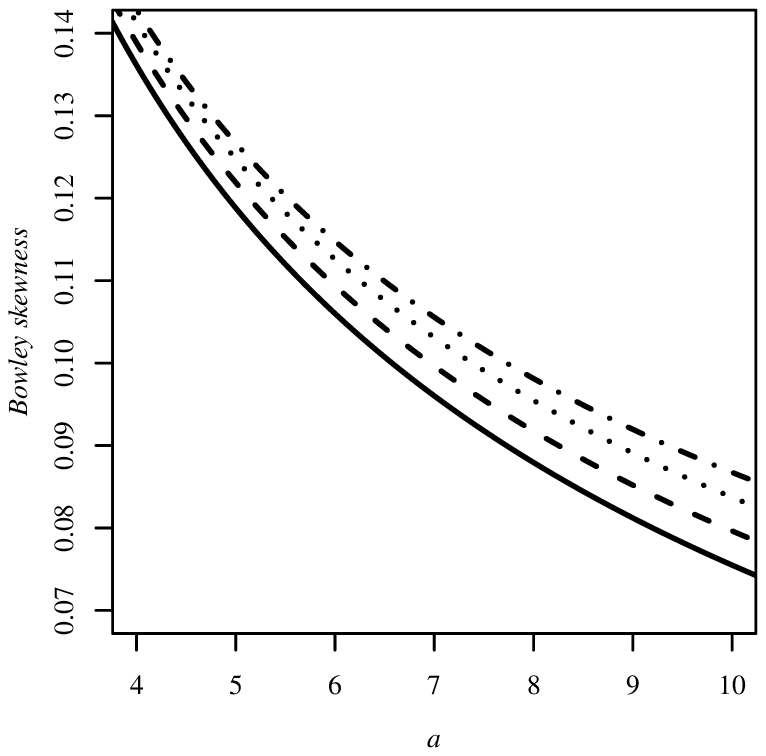,height=0.45\linewidth,angle=-90}}
\subfigure[\label{Figure4c}]{\epsfig{file=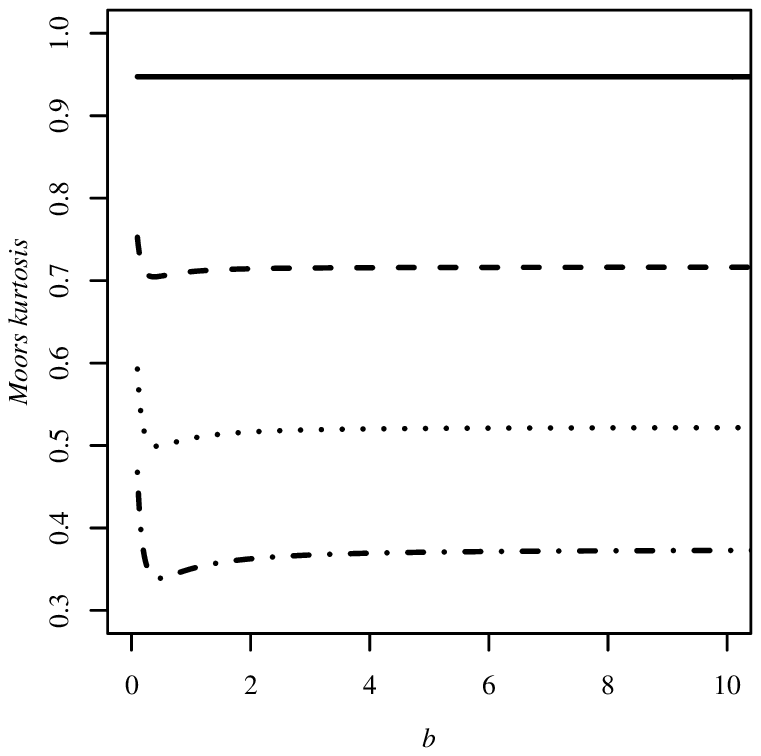,height=0.45\linewidth,angle=-90}}
\subfigure[\label{Figure4d}]{\epsfig{file=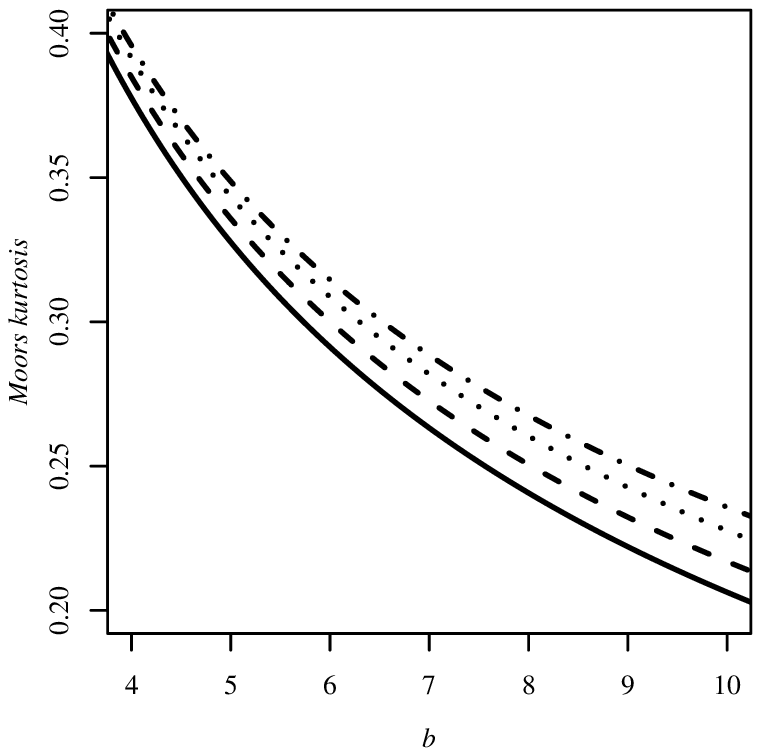,height=0.45\linewidth,angle=-90}}
\caption{Plots of the Bowley skewness and Moors kurtosis in terms of
(a)~$a$ for $b = 1.0$ (solid curve) and $b = 1.5$ (dashed curve),
$b = 3.5$ (dotted line) and $b = 4.5$ (bold line);
(b)~$b$ for $a =1.0$ (solid curve), $a = 1.5$ (dashed curve), $a = 3.5$ (dotted
line) and $a = 4.5$ (bold line);
(c)~$a$ for $b = 1.0$ (solid curve) and $b = 1.5$ (dashed curve), $b = 3.5$ (dotted line) and $b = 4.5$
(bold line);
and
(d)~$b$ for $a = 1.0$ (solid curve), $a = 1.5$
(dashed curve), $a = 3.5$ (dotted line) and $a = 4.5$ (bold line).}
\label{Figure4}
\end{figure}

Plots of the Bowley skewness and Moors kurtosis
for selected values of $a$ and $b$ are displayed in Figure~\ref{Figure4a}.
The parameter $\theta$ was set to one.

\section{Mean deviations and inequality measures}
\label{section-mean-deviations}

The amount of scatter in $X$ is measured to some extent by the totality of deviations from
the mean ($\mu$) and median ($m$). These are known as the mean deviation about the mean and
the mean deviation about the median given by
\begin{gather*}
\delta_1(X)
=
2\mu F(\mu) - 2\mu + 2\int_{\mu}^{\infty}x f(x)\mathrm{d}x
\quad
\text{and}
\quad
\delta_2(X)
=
2\int_{m}^{\infty}x f(x) \mathrm{d}x - \mu
,
\end{gather*}
respectively,
where
$\mu=\operatorname{E}(X)$ and
$m=Q(1/2)$.

Defining the integral $J(z)=\int^{z}_{0}xf(x) \mathrm{d}x$,
the measures
$\delta_1(X)$ and $\delta_2(X)$ are given by
\begin{align*}
\delta_1(X)=2\mu F(\mu)-2 J(\mu)
\quad
\text{and}
\quad
\delta_2(X)=\mu-2J(m)
,
\end{align*}
where $F(\mu)$ and $F(m)$ are easily obtained from~\eqref{eq2}.

We now determine $J(z)$.
Substituting $y=\theta x^{-2}$ in equation~\eqref{eq3},
we obtain
\begin{align*}
J(z)
=
\int_0^z x f(x) \mathrm{d}x =
\frac{\sqrt\theta}{B(a,b)}
\int_{\theta/z^2}^\infty
y^{-1/2}\exp(-ay)
\left[1 - \exp(-y)\right]^{b-1}
\mathrm{d}y
.
\end{align*}

Considering the power series~\eqref{eq7},
we have
\begin{align*}
J(z)
=&
\frac{\sqrt\theta}{B(a,b)}
\sum_{n=0}^\infty
\frac{(-1)^n}{\Gamma(b-n)n!}
\int_{\theta/z^2}^{\infty}
y^{-1/2}
\exp\{-(a+n)y \}
\mathrm{d}y
\\
=&
\sqrt{\pi\theta}
\frac{\Gamma(b)}{B(a,b)}
\sum_{n=0}^{\infty}
\frac{ (-1)^n}{\Gamma(b-n)n!}
\frac{1}{\sqrt{a+n}}
\operatorname{erfc}
(\sqrt{\theta/z^2}\sqrt{a+n}),
\end{align*}
where
$\operatorname{erfc}(x)=\frac{2}{\sqrt{\pi}}\int_{x}^{\infty}e^{-t^2}\mathrm{d}t$
is the complementary error function.

\cite{BoF:30} and \cite{Lor:05}
curves are inequality measures
which have applications
in economics, reliability, demography,
actuarial sciences, and medicine, among others.
They are defined by
\begin{align*}
B(p)
=
\frac{1}{p\mu}\int_{0}^{Q(p)}xf(x) \mathrm{d}x
=
\frac{1}{p\mu}J(Q(p))
\quad
\text{and}
\quad
L(p)
=
\frac{1}{\mu}\int_{0}^{Q(p)}xf(x) \mathrm{d}x
=
\frac{1}{\mu}J(Q(p))
,
\end{align*}
respectively, for $0<p\leq1$, see \cite{PuD:05} for details.

\section{Shannon and R\'enyi entropies}
\label{section-entropy}

The entropy of a random variable quantifies its associated uncertainty \citep{SoG:01}.
Two important entropy measures are the Shannon entropy and its generalization known as the
R\'enyi entropy. For the BIR distribution, the Shannon entropy is

\begin{align*}
H(X)
=
&
-\mathrm{E}\{\log[f(X)]\} =
-\int^{\infty}_{0}f(x) \log[f(x)]
\mathrm{d}x
\\
=
&
-
\int^{\infty}_{0}f(x)
\log \left\{
\frac{2\theta}{\operatorname{B}(a,b)\,x^{3}}
\exp\left(-\frac{a\theta}{x^2}\right)
\left[1-\exp\left(-\frac{\theta}{x^2}\right)\right]^{b-1}
\right\}
\mathrm{d}x
\\
=
&
-
\int^{\infty}_{0}
\log
\left(\frac{2\theta}{\operatorname{B}(a,b)}\right)
f(x)
\mathrm{d}x
+
3\int^{\infty}_{0} \log(x) f(x)
\mathrm{d}x
\\
&
+
a\theta\int^{\infty}_{0} \frac{1}{x^2} f(x)
\mathrm{d}x
-
(b-1)
\int^{\infty}_{0}
\log
\left[ 1-\exp\left(-\frac{\theta}{x^2}\right) \right]
f(x)
\mathrm{d}x
,
\end{align*}
where the first of the last four integrals is equal to $-\log
[2\theta/\operatorname{B}(a,b)]$. The second integral can be
calculated as follows

\begin{align*}
\int^{\infty}_{0}
\log(x) f(x)dx
=&
\int^{\infty}_{0}
\log(x)
\left\{
\frac{2\theta}{\operatorname{B}(a,b)\,x^{3}}
\exp\left(-\frac{a\theta}{x^2}\right)
\left[1-\exp\left(-\frac{\theta}{x^2}\right)\right]^{b-1}
\right\}
\mathrm{d}x
\\
=
&
\frac{2\theta}{\operatorname{B}(a,b)}
\sum_{n=0}^{\infty}
\frac{(-1)^n\Gamma(b)}{\Gamma(b-n)n!}
\int^{\infty}_{0}
\frac{\log(x)}{x^3}
\exp\left\{\frac{-(a+n)\theta}{x^2}\right\}
\mathrm{d}x
\\
=
&
\frac{2\theta}{\operatorname{B}(a,b)}
\sum_{n=0}^{\infty}
\frac{(-1)^n\Gamma(b)}{\Gamma(b-n)n!}
\frac{\log[(a+n)\theta]+\gamma}{4(a+n)\theta}
,
\end{align*}
where $\gamma$ is the Euler-Mascheroni constant. The third integral can be expressed by
\begin{align*}
\int^{\infty}_{0}
\frac{1}{x^2} f(x)
\mathrm{d}x
=
&
\int_{0}^{\infty}
\frac{2\theta}{\operatorname{B}(a,b)\,x^{5}}
\exp\left(-\frac{a\theta}{x^2}\right)
\left[1-\exp\left(-\frac{\theta}{x^2}\right)\right]^{b-1}
\mathrm{d}x
\\
=
&
\frac{2\theta}{\operatorname{B}(a,b)}
\int^{\infty}_{0}
\frac{1}{x^5}
\exp\left(-\frac{a\theta}{x^2}\right)
\sum_{n=0}^{\infty}
\frac{(-1)^n\Gamma(b)}{\Gamma(b-n)n!}
\exp\left\{\frac{-n\theta}{x^2}\right\}
\mathrm{d}x
\\
=
&
\frac{2\theta}{\operatorname{B}(a,b)}
\sum_{n=0}^{\infty}
\frac{(-1)^n\Gamma(b)}{\Gamma(b-n)n!}
\int^{\infty}_{0}
\frac{1}{x^5}
\exp\left\{\frac{-(a+n)\theta}{x^2}\right\}
\mathrm{d}x
\end{align*}

Setting
$t = \frac{(a+n)\theta}{x^2}$,
we obtain
\begin{align*}
\int^{\infty}_{0}
\frac{1}{x^2} f(x)
\mathrm{d}x
=
&
\frac{2\theta}{\operatorname{B}(a,b)}
\sum_{n=0}^{\infty}
\frac{(-1)^n\Gamma(b)}{\Gamma(b-n)n!}
\int^{\infty}_{0}
\frac{\exp(-t)}{[2(a+n)\theta]^2}
\mathrm{d}t
\\
=
&
\frac{1}{\theta\operatorname{B}(a,b)}
\sum_{n=0}^{\infty}
\frac{(-1)^n\Gamma(b)}{\Gamma(b-n)n!}
\frac{1}{(a+n)^2}
.
\end{align*}

Considering the fourth integral, let $u=\theta x^{-2}$. From the power series expansion
$\log(1+z)=z+\frac{1}{2}z^2- \frac{1}{3} z^3-\cdots$,
we can write
\begin{align*}
\int^{\infty}_{0}
\log[ 1-\exp(-\theta/x^2) ]f(x)
\mathrm{d}x
=
&
\frac{1}{\operatorname{B}(a,b)}
\int^{\infty}_{0}
\log\{1-\exp(-u)\}
\exp(-au)[1-\exp(-u)]^{b-1}
\mathrm{d}u
\\
=
&
-\frac{1}{\operatorname{B}(a,b)}
\int_{0}^{\infty}
\sum_{k=1}^{\infty}
\frac{\exp\{-u(k+a)\}
\left[1-\exp(-u)\right]^{b-1}}{k}
\mathrm{d}u
\\
=
&
\frac{1}{\operatorname{B}(a,b)}
\sum_{k=1}^{\infty}
\sum_{n=0}^{\infty}
\frac{(-1)^{n+1}\,\Gamma(b)}{k\Gamma(b-n)\,n!}
\int_{0}^{\infty}
\exp\left\{-u(a+k+n)\right\}
\mathrm{d}u
\\
=
&
\frac{1}{\operatorname{B}(a,b)}
\sum_{k=1}^{\infty}
\sum_{n=0}^{\infty}
\frac{(-1)^{n+1}\,\Gamma(b)}{k(a+k+n)\Gamma(b-n)\,n!}
.
\end{align*}

Finally, we obtain
\begin{equation*}
\begin{split}
H(X)
=
&
-\log
\left\{
\frac{2\theta}{\operatorname{B}(a,b)}
\right\}
+
\\
&
\frac{\Gamma(b)}{\operatorname{B}(a,b)}
\sum_{n=0}^{\infty}
\frac{(-1)^n}{\Gamma(b-n)\,n!}
\left[
\frac{3}{2}\frac{\log\{(a+n)\theta\}+\gamma}{a+n}
+
\frac{a}{(a+n)^2}
+
(b-1)
\sum_{k=1}^{\infty}
\frac{1}{k(a+k+n)}
\right]
.
\end{split}
\end{equation*}

Now, the R\'enyi entropy can be expressed as
\begin{align}
\label{equation-renyi}
H_\alpha(X)
=
\frac{1}{1-\alpha}
\log
\left(
\int^\infty_0
f(x)^{\alpha}
\mathrm{d}x
\right)
,
\quad
\alpha > 0,
\;
\alpha \neq 1,
\end{align}
where $\alpha > 0$ and $\alpha \neq 1$.
Notice that when $\alpha\to1$, the R\'enyi entropy converges to the Shannon entropy. For
calculating~\eqref{equation-renyi}, we apply~\eqref{eq3} and consider the power series
expansion~\eqref{eq7} yielding
\begin{equation*}
\label{enRny}
\begin{split}
\int^\infty_0
f(x)^\alpha
\mathrm{d}x
=
&
\left[\frac{2\theta}{\operatorname{B}(a,b)}\right]^\alpha
\Gamma(\alpha(b-1)+1)
\sum_{n=0}^\infty
\frac{(-1)^n}{\Gamma(\alpha(b-1)+1-n)\,n!}
\\
&
\times
\int^\infty_0
x^{-3\alpha}
\exp\left\{-(a\alpha+n)\frac{\theta}{x^2} \right\}
\mathrm{d}x
.
\end{split}
\end{equation*}

The last integral
can be evaluated as follows.
Let $u=\theta x^{-2}$. Then, we have

\begin{align*}
\int^\infty_0
x^{-3\alpha}
\exp\left\{-(a\alpha+n)\frac{\theta}{x^2} \right\}
\mathrm{d}x
=
&
\int^\infty_0
u^{\frac{3(\alpha-1)}{2}}
\exp\left\{-(a\alpha+j)u\right\}
\mathrm{d}u
\\
=
&
\left(
\frac{1}{a\alpha+n}
\right)^{\frac{3\alpha-1}{2}}
\Gamma\left( \frac{3\alpha-1}{2} \right)
.
\end{align*}

Finally, we obtain
\begin{align*}
\int^\infty_0 f(x)^\alpha \mathrm{d}x
=
\left[\frac{2\theta}{\operatorname{B}(a,b)}\right]^\alpha
\frac{\Gamma\left( \frac{3\alpha-1}{2}\right)}
{2\theta^{\frac{3(\alpha-1)}{2}+1}}
\sum^\infty_{n=0}
\frac{(-1)^n}{(a\alpha+n)^{\frac{3\alpha-1}{2}}}
\frac{\Gamma(\alpha(b-1)+1)}{\Gamma(\alpha(b-1)+1-n)n!}
.
\end{align*}

\section{Order statistics}
\label{section-order}

Here, we present an explicit expression for
the density function $f_{i:n}(x)$ of
the $i$th order statistic $X_{i:n}$
in a random sample of size $n$ from the
BIR distribution.
Consider the well-known result
\begin{align*}
f_{i:n}(x)=\frac{f(x)}{\operatorname{B}(i,n-i+1)}
F(x)^{i-1}\{1-F(x)\}^{n-i},
\end{align*}
for $i=1,\ldots,n$.
Applying the binomial expansion in the above equation,
we obtain
\begin{align*}
f_{i:n}(x)=
\frac{f(x)}{\operatorname{B}(i,n-i+1)}
\sum_{l=0}^{n-i}
\binom{n-i}{l}
(-1)^l
F(x)^{i+l-1}
.
\end{align*}

Inserting~\eqref{eq3} and \eqref{cdf-expansion} in the last equation,
$f_{i:n}(x)$ can be expressed as
\begin{equation}
\label{ord-04}
\begin{split}
f_{i:n}(x)
=&
\frac{2\theta}{B(i,n-i+1)x^3}
\exp\left(-\frac{a\theta}{x^2}\right)
\left[1-\exp\left(-\frac{\theta}{x^2}\right)\right]^{b-1}
\sum_{l=0}^{n-i}
\binom{n-i}{l}
\frac{(-1)^l}{\operatorname{B}(a,b)^{i+l}}
\\
&
\times
\left[\sum^{\infty}_{j=0} \frac{(-1)^j}{a+j}
\frac{\Gamma(b)}{\Gamma(b-j)j!}
\exp\left(-\frac{(a+j)\theta}{x^2}\right) \right]^{i+l-1},
\end{split}
\end{equation}
for $b>0$ real non-integer.

Now, using the following identity
\begin{eqnarray*}
\left(
\sum_{i=0}^\infty
a_{i}
\right)^{k}
=
\sum_{m_{1}=0}^{\infty}\cdots\sum_{m_{k}=0}^{\infty}a_{m_{1}}\cdots{a}_{m_{k}},
\end{eqnarray*}
for $k$ positive integer, we can write~\eqref{ord-04} as
\begin{equation}\label{ord-05}
 f_{i:n}(x)=\sum_{l=0}^{n-i}\sum_{m_{1}=0}^{\infty}\cdots\sum_{m_{i+l-1}=0}^{\infty}\delta_{i,l}f_{i,l}(x),
\end{equation}
where
\begin{eqnarray*}
f_{i,l}(x)
=
\frac{
2\theta\exp\left(-\frac{a\theta}{x^2}\right)
\left[1-\exp\left(-\frac{\theta}{x^2}\right)\right]^{b-1}
\exp
\left\{
-\frac{\theta}{x^2} \sum_{j=1}^{i+l-1}(a+m_j)
\right\}
}
{x^{3}\operatorname{B}\left( a(i+l)+\sum_{j=1}^{i+l-1}m_{j},b \right)},
\end{eqnarray*}
and
\begin{eqnarray*}
\delta_{i,l}
=
\frac{(-1)^{l+\sum_{j=1}^{i+l-1}m_{j}}{\binom{n-i}{l}}
\Gamma(b)^{i+l-1}
\operatorname{B}
\left(
a(i+l)+\sum_{j=1}^{i+l-1}m_{j},b
\right)}
      {{\operatorname{B}(a,b)^{i+l}} {\operatorname{B}(i,n-i+1)}\prod_{j=1}^{i+l-1}(a+m_{j})\Gamma(b-m_{j})m_{j}!}.
\end{eqnarray*}

Note that $f_{i,l}(x)$ is the density function of
the $\operatorname{BIR}(a(i+l)+\sum_{j=1}^{i+l-1},b,\theta)$ distribution.
Also, the constants $\delta_{i,l}$ are obtained given
$i,n,l$ and a sequence of indices $m_{1},\ldots,m_{i+l-1}$.
The sums in~\eqref{ord-05} extend over all
$(i+l)$-tuples $(l,m_{1},\ldots,m_{i+l-1})$ of non-negative
integers.
These sums indicate that the density function of the BIR order statistics is
a linear combination of BIR densities. So, several structural quantities of
the BIR order statistics can be obtained from those of BIR distribution.

\section{Maximum likelihood estimation and information matrix}
\label{section-mle}

Consider independent BIR distributed random variables
$X_{1},\ldots,X_{n}$
with parameter vector
$\bm{\lambda} = (a,b,\theta)^{T}$.
The log-likelihood function $\ell(\bm{\lambda})$
for the BIR model
reduces to
\begin{align*}
\ell(a,b,\theta)
=&
n
[
\log(2\theta) - \log \{ \operatorname{B}(a,b) \}
]
-
3 \sum^n_{i=1}
\log (x_{i})
-
a\theta
\sum^n_{i=1} \frac{1}{x^2_i}
\nonumber
\\
&
+
(b-1)
\sum^{n}_{i=1}
\log\left\{
1- \exp\left(-\frac{\theta}{x_i^2}\right)
\right\}
.
\end{align*}

The elements of the score vector are:
\begin{align*}
U_{a}(\bm{\lambda}) =&
\frac{\partial}{\partial a}
\ell(a,b,\theta)
=
n[\psi (a+b) -\psi (a) ]
-
\theta
\sum^n_{i=1}
\frac{1}{x_i^2}
,
\\
U_{b}(\bm{\lambda})
=&
\frac{\partial}{\partial b}
\ell(a,b,\theta)
=
n[\psi (a+b) - \psi (b)]
 +
\sum^n_{i=1}
\log\left\{
1- \exp\left(-\frac{1}{x_i^2}\right)
\right\}
,
\\
U_{\theta}(\bm{\lambda})
=&
\frac{\partial}{\partial \theta}
\ell(a,b,\theta)
=
\frac{n}{\theta}
-
\sum^n_{i=1}
\frac{1}{x_i^2}
\left[
a
+
\frac{b-1}
{1-   \exp\left(\frac{\theta}{x_i^2} \right) }
\right]
,
\end{align*}
where $\psi(\cdot)$ is the digamma function, see \cite{AbS:72}.
The maximum likelihood equations can be solved numerically
for $a$, $b$, and $\theta$.
Under standard regularity conditions \citep{Cox:74} that are fulfilled for the proposed model whenever the parameters are in the interior of the parameter space,  the observed information matrix~$\mathcal{I}(\bm \lambda)$
can be employed for interval estimation of the model parameters and for hypothesis tests. The
BIR observed information matrix is given by
\begin{align*}
\mathcal{I}(\bm{\lambda})=
-
\begin{bmatrix}
U_{aa}(\bm{\lambda})&U_{ab}(\bm{\lambda})&U_{a\theta}(\bm{\lambda})  \\
U_{ab}(\bm{\lambda})&U_{bb}(\bm{\lambda})&U_{b\theta}(\bm{\lambda})  \\
U_{a\theta}(\bm{\lambda})&U_{b\theta}(\bm{\lambda})&U_{\theta\theta}(\bm{\lambda})
\end{bmatrix}
,
\end{align*}
whose elements are
\begin{align*}
U_{aa}(\bm{\lambda})
=
\frac{\partial}{\partial a}
U_a(\bm{\lambda})
=&
n \left[\psi_1(a+b) - \psi_1(a) \right]
,
\\
U_{ab}(\bm{\lambda})
=
\frac{\partial}{\partial b}
U_a(\bm{\lambda})
=&
n \psi_1(a+b), \\
U_{a\theta}(\bm{\lambda})
=
\frac{\partial}{\partial \theta}
U_a(\bm{\lambda})
=&
-
\sum^n_{i=1}
\frac{1}{x_i^2}
,
\\
U_{bb}(\bm{\lambda})
=
\frac{\partial}{\partial b}
U_b(\bm{\lambda})
=&
n \left[\psi_1(a+b) - \psi_1(b) \right]
,
\\
U_{b\theta}(\bm{\lambda})
=
\frac{\partial}{\partial \theta}
U_b(\bm{\lambda})
=&
-\sum^n_{i=1}
\frac{1 }
{x^2_i  \left[1-\exp\left(\frac{\theta}{x_i^2} \right) \right]}
,
\\
U_{\theta\theta}(\bm{\lambda})
=
\frac{\partial}{\partial \theta}
U_\theta(\bm{\lambda})
=&
-\frac{n}{\theta^2}
+
(b-1)
\sum^{n}_{i=1}
\frac{1}{x_i^4}
\frac{1}
{\left[1-\exp\left(\frac{-\theta}{x_i^2} \right) \right]
\cdot
\left[1-\exp\left(\frac{\theta}{x_i^2} \right) \right]}
,
\end{align*}
and $\psi_1(\cdot)$ is the polygamma function,
which satisfies
$\psi_1(x) = \frac{\mathrm{d}}{\mathrm{d}x} \psi(x)$.

\section{Application to real data}
\label{section-application}

\cite{Bjr:60} considered guinea pigs as a model to study human tuberculosis.
Indeed, guinea pigs are highly susceptible to the associate pathogen:
tubercle bacillus.
Bjerkedal generated several data sets with the
survival times of infected guinea pigs.
A particular set with 72 observations is listed in Table \ref{Table1}.
We analyze these data and fit to them the BIR,
exponentiated inverse Rayleigh (EIR) \citep{Gut:98},
IR, Rayleigh (R) \citep{t:64} and generalized Rayleigh (GR) \citep{vd:76}; \citep{VDa:76}
distributions.

All these distributions are common models for lifetime data.
Computational implementation was performed in
\texttt{Ox} matrix programming language \citep{Doo:06}.

Table \ref{Table2} lists the maximum likelihood estimates (MLEs)
of the model parameters (standard errors in parentheses) for each model.
It is also shown the values for
the Akaike information criterion (AIC) \citep{Akk:73}.
Bayesian information criterion (BIC) \citep{SCH:78},
bias-corrected Akaike information criterion (BAIC) \citep{HTs:89}
and Hannan-Quinn information criterion (HQIC) \citep{HQn:79}.
These results indicate that the BIR distribution has the lowest AIC, BAIC, and HQIC values among
the current models. Moreover, the GR distribution presents the lowest BIC value.

\begin{table}[h!]
\centering
\caption{Survival times of guinea pigs injected with tubercle bacilli}
\bigskip
\begin{tabular}{@{}ccccccccccccccc} \hline %
 12  & 15  & 22  & 24  & 24  & 32  & 32  & 33  & 34  & 38  & 38  & 43  & 44  & 48  & 52 \\%\midrule
 53  & 54  & 54  & 55  & 56  & 57  & 58  & 58  & 59  & 60  & 60  & 60  & 60  & 61  & 62 \\
 63  & 65  & 65  & 67  & 68  & 70  & 70  & 72  & 73  & 75  & 76  & 76  & 81  & 83  & 84 \\
 85  & 87  & 91  & 95  & 96  & 98  & 99  & 109 & 110 & 121 & 127 & 129 & 131 & 143 & 146 \\
 146 & 175 & 175 & 211 & 233 & 258 & 258 & 263 & 297 & 341 & 341 & 376 &     &     &     \\ \hline%
\end{tabular}
\label{Table1}
\end{table}

\begin{table}[h]
\centering
\caption{MLEs of the Model Parameters for the Survival Time Data.}\label{Table2}
\bigskip
\begin{tabular}{c@{\quad }c@{\,\,\,}c@{\,\,\,}c  c@{\quad }c@{\,\,\,}c c@{\quad}c@{\,\,\,}c}\hline
         & \multicolumn{4}{c}{Estimates} & \multicolumn{5}{c}{Goodness-of-fit} \\
& \multicolumn{4}{c}{(standard errors)} & \multicolumn{5}{c}{measures} \\ \hline

& $\widehat{a}$  & $\widehat{b}$              & $\widehat{\theta}$ & $\widehat{\lambda}$   & &  AIC  &BIC &BAIC &HQIC  \\ \hline

BIR               &\multicolumn{1}{r}{1094.47} & 0.61666             & 1.23294               & {$-$} & & {805.83}  & {812.66} & {878.19} & {808.55} \\
&\multicolumn{1}{r}{(349.183)} &(0.08653)  &(0.30105)  & {$(-)$} \\

EIR               &\multicolumn{1}{r}{101.482} & 1.00000             & 21.5592               & {$-$} & & {817.47}  & {822.03} & {889.65} &  {819.28} \\
&\multicolumn{1}{r}{(155.115)} &$(-)$  &(0.15736)  & {$(-)$} \\

IR               &{$-$} & $-$            & 2187.88               &{$-$}                              & & {815.47}   & {817.75} & {887.53} & {816.38} \\
&        {$(-)$} &$(-)$  &(257.844)  &{$(-)$} \\

R                   &{$-$} & $-$            & 90.6963               &  {$-$}                         & &{818.59}  & {820.87} & {890.65} & {819.50} \\
&           {$(-)$} &$(-)$  &(969.36)  &{$(-)$} \\

GR               &{$-$} & $-$            & 0.59904               &\multicolumn{1}{c}{0.00645}        & & {807.04}  & {811.60} & {879.22} & {808.86} \\
&        {$(-)$} &$(-)$  &(0.08742)  &\multicolumn{1}{r}{(0.00059)} \\ \hline
\end{tabular}
\end{table}

\begin{figure}[t]
\centering
\epsfig{file=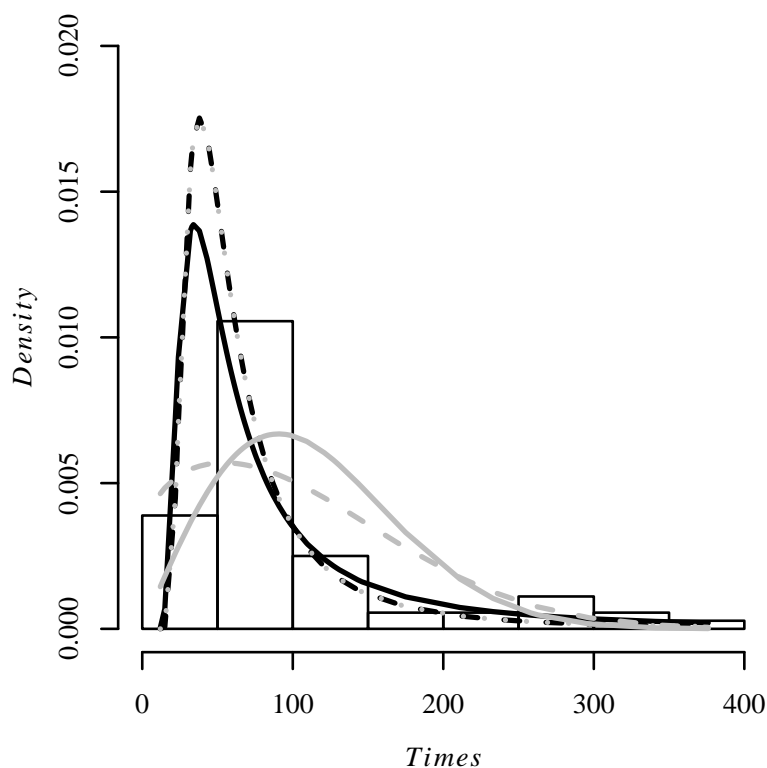,height=0.65\linewidth,angle=-90}
\caption{%
Fitted densities for
the BIR (black solid),
EIR (black dashed),
IR (gray dotted),
Rayleigh (gray solid), and
GR (gray dashed)
models for the survival time data.}
\label{Figure5}
\end{figure}

Plots of the estimated densities of the BIR, EIR, IR, Rayleigh and
GR models fitted to these data are displayed in Figure $\ref{Figure5}$.
The overall results suggest that the BIR distribution is superior to the
remaining distributions in terms of model fitting.

\section{Conclusion}
\label{section-conclusion}

In this work, we study the beta inverse Rayleigh distribution as a generalization of the inverse
Rayleigh distribution. We also provide a better foundation for some mathematical  properties for this distribution, including the
derivation of the hazard rate function, moments, quantile measures, mean deviations, entropy measures
and order statistics.
The model parameters are estimated by maximum likelihood.
An application of the BIR distribution to a real data set indicates
that the new distribution outperforms several distributions,
including the IR and Rayleigh distributions.

\section*{Acknowledgements}
The authors acknowledge support from CAPES, CNPq, and FACEPE.

\end{document}